\def\marker{\>\hbox{${\vcenter{\vbox{
    \hrule height 0.4pt\hbox{\vrule width 0.4pt height 6pt
    \kern6pt\vrule width 0.4pt}\hrule height 0.4pt}}}$}\>}
\def\gpic#1{#1
     \medskip\par\noindent{\centerline{\box\graph}} \medskip}
\documentclass{article}
\usepackage{jeffe}
\usepackage{url}
\usepackage{verbatim}
\usepackage{subfigure}
\usepackage{graphics}
\usepackage{fullpage}
\newtheorem{theorem}{Theorem}
\newtheorem{conjecture}[theorem]{Conjecture}
\newtheorem{lemma}[theorem]{Lemma}

\def\floor#1{\lfloor #1 \rfloor}

\begin{document}

\title{Edge-choosability and total-choosability of \\planar graphs with no adjacent 3-cycles}
\author{Daniel W. Cranston\footnote{DIMACS, Rutgers University, Piscataway, NJ} \\ \texttt{dcransto@dimacs.rutgers.edu}}
\maketitle
\begin{abstract}
Let $G$ be a planar graph with no two 3-cycles sharing an edge. 
We show that if $\Delta(G)\geq 9$, then $\chi'_l(G) = \Delta(G)$ and $\chi''_l(G)=\Delta(G)+1.$
We also show that if $\Delta(G)\geq 6$, then $\chi'_l(G)\leq\Delta(G)+1$  
and if $\Delta(G)\geq 7$, then $\chi''_l(G)\leq\Delta(G)+2$. 
All of these results extend to graphs in the projective plane and when $\Delta(G)\geq 7$ the results also extend to graphs in the torus and Klein bottle.

This second edge-choosability result improves on work of Wang and Lih and of Zhang and Wu.
All of our results use the discharging method to prove structural lemmas about the existence of subgraphs with small degree-sum.
For example, we prove that if $G$ is a planar graph with no two 3-cycles sharing an edge and with $\Delta(G)\geq 7$, 
then $G$ has an edge $uv$ with $d(u)\leq 4$ and $d(u)+d(v)\leq \Delta(G)+2$.

All of our proofs yield linear-time algorithms that produce the desired colorings.
\bigskip

\noindent
\textit{MSC}: 05C15, 05C10
\bigskip

\noindent
\textbf{Keywords}: List coloring, edge coloring, total coloring, Vizing's Conjecture
\end{abstract}

\smallskip

\section{Introduction}
All our graphs are finite and without loops or multiple edges.  Let $G$ be a plane graph.  
We use $E(G)$, $V(G)$, $F(G)$, $\Delta(G)$, and $\delta(G)$ to denote the edge set, vertex set, face set, maximum degree, and minimum degree of $G$, respectively.  
When the graph is clear from context, we use $\Delta$, rather than $\Delta(G)$.
We use ``$j$-face'' and ``$j$-vertex'' to mean faces and vertices of degree $j$.  The degree of a face $f$ is the number of edges along the boundary of $f$, with each cut-edge being counted twice.  The degree of a face $f$ and the degree of a vertex $v$ are denoted by $d(f)$ and $d(v)$.  
We say a face $f$ or vertex $v$ is \textit{large} when $d(f)\geq 5$ or $d(v)\geq 5$.
We use \textit{triangle} to mean 3-cycle.
We use \textit{kite} to mean a subgraph of $G$ formed by two 3-cycles that share an edge.
We use \textit{element} to mean vertex or face.

A \textit{proper total-coloring} of $G$ is an assignment of a label to each element so that no two incident or adjacent elements receive the same label.  
We call these labels \textit{colors}.  
A \textit{proper $k$-total-coloring} is a proper total-coloring that uses no more than $k$ colors.
A \textit{total assignment} $L$ is a function on $E(G)\cup V(G)$ that assigns each element $x$ a list $L(x)$ of colors available for use on that element.  
An \textit{$L$-total-coloring} is a proper total-coloring with the additional constraint that each element receives a color appearing in its assigned list.  
We say that a graph $G$ is \textit{$k$-total-choosable} if $G$ has a proper $L$-total-coloring whenever $|L(x)|\geq k$ for every $x\in E(G)\cup V(G)$.
The \textit{total chromatic number} of $G$, denoted $\chi''(G)$, is the least integer $k$ such that $G$ is $k$-total-colorable.  
The \textit{list total chromatic number} of $G$, denoted $\chi''_l(G)$, is the least integer $k$ such that $G$ is $k$-total-choosable.  
In particular, note that $\chi''(G)\leq \chi''_l(G)$.  
The \textit{list edge chromatic number} $\chi'_l(G)$ is defined similarly in terms of coloring only edges; the ordinary edge chromatic number is denoted $\chi'(G)$.
Probably the most fundamental and important result about the edge chromatic number of graphs is:

\begin{theorem}
(Vizing's Theorem; Vizing \cite{vizing1, vizing2} and Gupta \cite{gupta})
\label{vizing-thm}
$$
\chi'(G)\leq\Delta(G)+1.
$$
\end{theorem}

Vizing conjectured that Theorem \ref{vizing-thm} could be strengthened by proving the same bound for the list edge chromatic number:
\begin{conjecture} 
(Vizing's Conjecture; see \cite{kostochka})
\label{vizing-conj}
$$
\chi'_l(G)\leq\Delta(G)+1.
$$
\end{conjecture}

The most famous open problem about list edge-coloring is the List Coloring Conjecture.  
Bollob\'{a}s and Harris \cite{bollobas} believed that Vizing's conjecture could be further strengthened to give:

\begin{conjecture}
(List Coloring Conjecture; Bollob\'{a}s and Harris \cite{bollobas})
\label{LCC}
$$
\chi'_l(G) = \chi'(G).
$$
\end{conjecture}

We give a brief summary of previous results on list edge-coloring; for a more thorough treatment, we recommend \textit{Graph Coloring Problems}~\cite{JT}.
Borodin, Kostochka, and Woodall \cite{BKW} proved that the List Coloring Conjecture holds for planar graphs with $\Delta\geq 12$.  
Vizing's Conjecture is easy to prove when $\Delta\leq 2$.  
Harris \cite{harris} and Juvan et al.\ \cite{juvan1} confirmed the conjecture when $\Delta=3$ and $\Delta = 4$, respectively.  
Borodin proved Vizing's Conjecture for planar graphs with $\Delta\geq 9$~\cite{borodin}.  
Wang and Lih \cite{wang2} proved that Vizing's Conjecture holds for a planar graph $G$ when $\Delta\geq 6$ and $G$ has no two triangles sharing a vertex.  
Zhang and Wu \cite{zhang} proved that Vizing's Conjecture holds for a planar graph $G$ when $\Delta\geq 6$ and $G$ has no 4-cycles.  

We improve these results in several ways.  In Section 2, we prove structural results for use in Section 4 and Section 5, where we prove our main results.
For simplicity, we state each of our results only for planar graphs.  However, in Section 3 we show that each result can be extended to the projective plane and that most of the results can also be extended to the torus and Klein bottle.
In Section 4, we show that Vizing's Conjecture holds for a planar graph that contains no kites and has $\Delta\geq 6$.  This is a strengthening of the result of Wang and Lih \cite{wang2} and the result of Zhang and Wu \cite{zhang}.  
We also show that the List Coloring Conjecture holds for a planar graph that contains no kites and has $\Delta\geq 9$.
%
In Section 5 we prove results about list total coloring, which we describe below.

Less is known about the total chromatic number than the edge chromatic number.
Vizing and Behzad conjectured an analogue to Vizing's Theorem:

\begin{conjecture}
(Total Coloring Conjecture; Vizing~\cite{vizing1} and Behzad \cite{behzad})
$$
\chi''(G) \leq \Delta(G)+2.
$$
\end{conjecture}

The Total Coloring Conjecture was proved for $\Delta=3$ by Rosenfeld~\cite{rosenfeld} and also by Vijayaditya~\cite{vijayaditya}.  For $\Delta=4$ and $\Delta=5$ it was proved by Kostochka~\cite{kostochka2, kostochka3, kostochka4}.  
For planar graphs, much more is known.  Borodin~\cite{borodin2} proved the Total Coloring Conjecture for $\Delta\geq 9$.  
Yap~\cite{yap,JT} observed that the cases of the Total Coloring Conjecture when $\Delta=7$ or $\Delta=8$ follow from a short argument that uses the 4-Color Theorem and the fact that $\chi'(G)=\Delta$ for planar graphs when $\Delta\geq7$.
Borodin, Kostochka, and Woodall~\cite{BKW} showed that $\chi''(G)=\Delta+1$ for $\Delta\geq 12$.  

The list total chromatic number seems to have been relatively unstudied until Borodin, Kostochka, and Woodall conjectured the following:

\begin{conjecture}
(Total List Coloring Conjecture; Borodin, Kostochka, Woodall~\cite{BKW})
$$
\chi''_l(G) =\chi''(G).
$$
\end{conjecture}
For a planar graph with $\Delta(G)\geq 12$, they showed the stronger result $\chi''_l(G)=\chi''(G)=\Delta+1$.  We note that Borodin's proof of the Total Coloring Conjecture for planar graphs with $\Delta(G)\geq 9$ in fact shows that $\chi''_l(G)\leq\Delta+2$.

Almost all of our proofs for $\chi'_l(G)$ can easily be adapted to give results for $\chi''_l(G)$.  Again, we consider planar graphs with no kites.  In Section~5, we show that if $\Delta\geq 9$, then $\chi''_l(G)=\Delta+1$. We also show that if $\Delta\geq 7$, then $\chi''_l(G)\leq\Delta+2$.  


\section{Structure of planar graphs with no triangles sharing an edge}
Before we prove structural results about planar graphs with no kites, we mention a theorem of O.V. Borodin.

\begin{theorem}(Borodin~\cite{borodin3})
\label{borodin2thm}
Let $G$ be a plane graph with no kites.  The following statements are valid and all numerical parameters are best possible.
\begin{enumerate}
\item[(i)] $\delta(G)\leq 4$;
\item[(ii)] If $\delta(G)\geq 3$, then there are adjacent vertices $x$, $y$ such that $d(x)+d(y)\leq 9$;
\item[(iii)] if $\delta(G)\geq 3$, then there is either an $i$-face where $4\leq i\leq 9$ or 10-face incident with ten 3-vertices and adjacent to five triangles.
\end{enumerate}
\end{theorem}
We will not directly use Theorem~\ref{borodin2thm}, although part (ii) is very similar to some of our lemmas.

\begin{theorem}
\label{delta7}
If graph $G$ is planar, $G$ contains no kites, and $G$ has $\Delta\geq~7$, then $G$ has an edge $uv$ with $d(u)\leq 4$ and $d(u)+d(v)\leq \Delta + 2$.
\end{theorem}
\begin{proof}
Assume $G$ is a counterexample.  Clearly, $\delta(G)\geq 3$; furthermore, each 3-vertex is only adjacent to $\Delta$-vertices.
We use a discharging argument.  
We assign to each element $x$ an initial charge $\mu(x) = d(x) - 4$.  The initial charges have sum $\sum_{x\in V\cup F}(d(x)-4)=-8$.
We use the following two discharging rules, applied simultaneously at all vertices and faces in a single discharging phase:
\begin{itemize}
\item[]
\begin{itemize}
\item[(R1)] Each large vertex gives a charge of $1/2$ to each incident triangle.
\item[(R2)] Each $\Delta$-vertex gives a charge of $1/3$ to each adjacent 3-vertex.
\end{itemize}
\end{itemize}

To reach a contradiciton, we show that for every element the new charge $\mu^*$ is nonnegative.

\noindent Consider an arbitrary face $f$.  
\begin{itemize}
\item If $d(f) = 3$, then at least two of the vertices incident to $f$ are large; otherwise, we have $d(u)+d(v)\leq 4 + 4\leq\Delta+2$.
Thus $\mu^*(f)\geq -1 + 2(1/2) = 0$.  
\item If $d(f)\geq 4$, then $\mu^*(f) = \mu(f)\geq 0$.  
\end{itemize}
Consider an arbitrary vertex $v$.
\begin{itemize}
\item If $d(v) = 3$, then $\mu^*(v) = -1 + 3(1/3) = 0$, since each neighbor of $v$ is a $\Delta$-vertex.
\item If $d(v) = 4$, then $\mu^*(v) = \mu(v) = 0$.  
\item If $d(v) = 5$, then $v$ is incident to at most 2 triangles, so $\mu^*(v) \geq 1 - 2(1/2) = 0$.  
\item If $6\leq d(v) \leq \Delta-1$, then $v$ is incident to at most $d(v)/2$ triangles.  Thus $\mu^*(v)\geq d(v) -4 - \frac{d(v)}{2}\frac{1}{2} = \frac{3d(v)}{4} - 4 > 0$.  
\item If $d(v) = \Delta$, then let $t$ be the number of triangles incident to $v$. 
For each triangle incident to $v$, at most one of the vertices of that triangle has degree 3.  
Thus, if $v$ is incident to $t$ triangles, then $\mu^*(v)\geq {d(v)-4 -t(\frac{1}{2}) - (d(v)-t)(\frac{1}{3})} = {\frac{2d(v)}{3}- 4 - \frac{t}{6}}$.  
Since $t\leq \frac{d(v)}{2}$, we get $\mu^*(v) \geq \frac{7d(v)}{12} - 4$.  This expression is positive when $d(v)\geq 7$.
\end{itemize}
\aftermath
\end{proof}

We will use Theorem \ref{delta7} to show that any planar graph with $\Delta\geq 7$ that contains no kites is $(\Delta+1)$-edge-choosable.  We would also like to prove an analogous result for the case $\Delta=6$.  
To prove such a result, we need the following structural lemma.
We say that a triangle is of \textit{type} $(a,b,c)$ if its vertices have degrees $a, b,$ and $c$. 
%
\begin{lemma}
\label{delta6}
If graph $G$ is planar, $G$ contains no kites, and $\Delta = 6$,
then at least one of the three following conditions holds:
\begin{enumerate}
\item[(i)] $G$ has an edge $uv$ with $d(u)+d(v)\leq 8$.
\item[(ii)] $G$ has a 4-face $uvwx$ with $d(u) = d(w) = 3$.
\item[(iii)] $G$ has a 6-vertex incident to three triangles; two of these triangles are of type $(6,6,3)$ and the third is of type $(6,6,3)$, $(6,5,4)$, or $(6,6,4)$.
\end{enumerate}
\end{lemma}
\begin{proof}
Assume $G$ is a counterexample.  For every edge $uv$, $G$ must have $d(u) + d(v)\geq 9$.  Thus, $\delta(G)\geq 3$.
We use a discharging argument.  
We assign to each element $x$ an initial charge $\mu(x) = d(x) - 4$.
We use the following three discharging rules:
\begin{itemize}
\item []
\begin{itemize}
\item[(R1)] Each large face $f$ gives a charge of $1/2$ to each incident 3-vertex.
\item[(R2)] Each 5-vertex $v$ gives a charge of $1/2$ to each incident triangle.
\item[(R3)] Each 6-vertex $v$
\begin{itemize}
\item[$\bullet$] gives a charge of $1/3$ to each adjacent 3-vertex that is not incident to any large face.
\item[$\bullet$] gives a charge of $1/6$ to each adjacent 3-vertex that is incident to a large face.
\item[$\bullet$] gives a charge of $1/2$ to each incident triangle that is incident to a 3-vertex or a 4-vertex.
\item[$\bullet$] gives a charge of $1/3$ to each incident triangle that is not incident to a 3-vertex or a 4-vertex.
\end{itemize}
\end{itemize}
\end{itemize}

Now we show that for every element the new charge $\mu^*$ is nonnegative.

\noindent Consider an arbitrary face $f$.  
\begin{itemize}
\item If $d(f) = 3$, then we consider two cases.  If $f$ is incident to a 3-vertex or a 4-vertex, then $\mu^*(f) = -1 + 2(1/2) = 0$.  If $f$ is not incident to a 3-vertex or a 4-vertex, then $\mu^*(f) \geq -1 + 3(1/3) = 0$.  
\item If $d(f) = 4$, then $\mu^*(f) = \mu(f) = 0$.  
\item If $d(f) = 5$, then $\mu^*(f) \geq 1 - 2(1/2) = 0$.  
\item If $d(f)\geq 6$, then $\mu^*(f) \geq d(f) -4 - \frac{d(f)}{2}\frac{1}{2}  = \frac{3d(f)}{4} - 4 > 0$.
\end{itemize}

\noindent Consider an arbitrary vertex $v$.  
\begin{itemize}
\item If $d(v) = 3$, then we consider two cases.  If $v$ is incident to a large face, then $\mu^*(v) \geq -1 + 1/2 + 3(1/6) = 0$.  If $v$ is not incident to a large face, then $\mu^*(v) = -1 + 3(1/3) = 0$.  
\item If $d(v) = 4$, then $\mu^*(v) = \mu(v) = 0$.  
\item If $d(v) = 5$, then $\mu^*(v) \geq 1 - 2(1/2) = 0$.
\item 
If $d(v) = 6$, then we consider separately the four cases where $v$ is incident to zero, one, two, or three triangles.  
Note that if $v$ is incident to $t$ triangles, then the number of 3-vertices adjacent to $v$ is at most $(6-t)$.
\begin{itemize}
\item[$\circ$] If $v$ is incident to no triangles, then $\mu^*(v) \geq 2 - 6(1/3) = 0$.  
\item[$\circ$] If $v$ is incident to one triangle, then we consider two cases.
If $v$ is adjacent to at most four 3-vertices, then $\mu^*(v)\geq 2- (1/2) - 4(1/3)>0$.  If $v$ is adjacent to five 3-vertices, then two of these adjacent 3-vertices lie on a common face, together with $v$.  Since condition (ii) of the present lemma does not hold, this face must be a large face.  So $\mu^*(v) \geq 2 - (1/2) - 3(1/3) - 2(1/6) > 0$.  
\item[$\circ$] If $v$ is incident to two triangles, then we consider two cases.  If $v$ is adjacent to at most three 3-vertices, then $\mu^*(v)\geq 2 - 2(1/2) - 3(1/3) = 0$.  If $v$ is adjacent to four 3-vertices, then two of these adjacent 3-vertices lie on a common face, together with $v$.  Since condition (ii) of the present lemma does not hold, this face must be a large face.  So $\mu^*(v) \geq 2 - 2(1/2) - 2(1/3) - 2(1/6) = 0$.  
\item[$\circ$] If $v$ is incident to three triangles, then we consider two cases. If at most one of the triangles is type $(6,6,3)$, then $\mu^*(v) \geq 2 - 3(1/2) - 1/3 > 0$.  Furthermore, if two of the triangles incident to $v$ are type $(6,6,3)$ but the third triangle is not incident to any vertex of degree at most 4, then $\mu^*(v) = 2 - 2(1/2) - 2(1/3) - 1(1/3)~=~0$.  If two of the triangles are of type $(6,6,3)$ and the third triangle is incident to a vertex of degree at most 4, then condition (iii) of the lemma holds.
\end{itemize}
\end{itemize}
\vspace{-.1 in}
\end{proof}

We will apply Theorem~\ref{delta7}~and Lemma~\ref{delta6} to get our first result about edge-choosability.
To prove the $(\Delta+1)$-edge-choosability of a planar graph $G$ that has $\Delta\geq 6$ and that contains no kites, we remove one or more edges of $G$, inductively color the resulting subgraph, then extend the coloring to $G$.  
Intuitively, Theorem~\ref{delta7} and Lemma~\ref{delta6} do the ``hard work.''  
However, it is still convenient to prove the following lemma, which we will apply to the subgraphs of $G$ that arise from this process.

\begin{lemma}
\label{delta5-6}
Let $G$ be a planar graph that contains no kites. If $\Delta\leq5$, then $G$ has an edge $uv$ with ${d(u) + d(v)\leq 8}$.  If $\Delta=6$, then $G$ has an edge $uv$ with ${d(u) + d(v)\leq 9}$.
\end{lemma}
\begin{proof}
If $\Delta\leq 4$, then each edge $uv$ satisfies $d(u) + d(v)\leq 2\Delta \leq 8$.  
In that case, the lemma holds trivially.
If $\Delta=6$, then the result follows from $(ii)$ of Theorem~\ref{borodin2thm}.
So we must prove the lemma for the case $\Delta = 5$.
We use a discharging argument.
Assume $G$ is a counterexample.
For every edge $uv$, $G$ must have $d(u)+d(v)\geq \Delta + 4$.  Thus, $\delta(G)\geq 4$.  We assign to each element $x$ an initial charge $\mu(x) = d(x) - 4$.  We use a single discharging rule:
\begin{itemize}
\item[]
\begin{itemize}
\item[(R1)] Every large vertex $v$ gives a charge of $1/2$ to each incident triangle.
\end{itemize}
\end{itemize}

Now we show that for every element the new charge $\mu^*$ is nonnegative.

\noindent Consider an arbitrary face $f$.

\begin{itemize}
\item If $d(f)=3$, then $f$ is incident to at least two large vertices, so $\mu^*(f)\geq -1 + 2(1/2) = 0$.  
\item If $d(f)\geq 4$, then $\mu^*(f) = \mu(f) \geq 0$.  
\end{itemize}

\noindent Consider an arbitrary vertex $v$.
\begin{itemize}
\item If $d(v) = 4$, then $\mu^*(v) = \mu(v) =0$.  
\item If $d(v) = 5$, then $\mu^*(v)\geq 1 - 2(1/2) = 0$.  
\end{itemize}
\aftermath
\end{proof}
Note that the case when $\Delta=6$ in Lemma~\ref{delta5-6} follows easily from (ii) in Theorem~\ref{borodin2thm}; however, we reproved it above because in Section~3 we will adapt the proof of Lemma~\ref{delta5-6} to the projective plane, torus, and Klein bottle.

Before we state the next theorem, we need a new definition.  A $\textit{k-alternating cycle}$ is an even cycle $v_1w_1v_2w_2\ldots v_lw_l$ with $d(w_i)=k$.  This definition was introduced by Borodin in 1989~\cite{borodin2}.  Since then, it has been used to prove many coloring results (for example,~\cite{borodin}).

\begin{theorem}
\label{delta9}
If graph $G$ is planar, $G$ contains no kites, and $\Delta\geq 9$, then at
least one of the following two conditions holds:
\begin{enumerate}
\item[(i)] $G$ has an edge $uv$ with $d(u)\leq 4$ and $d(u)+d(v)\leq \Delta+1$.
\item[(ii)] $G$ has a $2$-alternating cycle $v_1w_1v_2w_2\ldots v_kw_k$. 
\end{enumerate}
\end{theorem}
\begin{proof}
Assume $G$ is a counterexample.  Clearly, $\delta(G)\geq 2$.  
Our proof will use a discharging argument, but first we show that 
if $G$ is a counterexample to Theorem~\ref{delta9}, then $G$ has more $\Delta$-vertices than 2-vertices.

Let $H$ be the subgraph of $G$ formed by all edges with one endpoint of degree 2 and the other endpoint of degree $\Delta$.  
Form $\widehat{H}$ from $H$ by contracting one of the two edges incident to each vertex of degree 2 (recall that each neighbor
of a 2-vertex in $G$ is a $\Delta$-vertex).  
Each 2-vertex in $G$ corresponds to an edge in $\widehat{H}$ and each vertex in $\widehat{H}$ corresponds to a $\Delta$-vertex in $G$.
So $G$ has more $\Delta$-vertices than 2-vertices unless $|E(\widehat{H})|\geq |V(\widehat{H})|$.  

If $|E(\widehat{H})|\geq |V(\widehat{H})|$, then $\widehat{H}$ contains a cycle.  
However, a cycle in $\widehat{H}$ corresponds to a 2-alternating cycle in $G$. 
Such a cycle in $G$ satisfies condition (ii) and shows that $G$ is not a counterexample to
Theorem~\ref{delta9}.  So, $G$ has more $\Delta$-vertices than 2-vertices.

We assign to each element $x$ an initial charge $\mu(x) = d(x) - 4$.
In addition to the vertices and edges, we create a \textit{bank} that can give and receive charge.
The bank has initial charge 0.  As with the vertices and edges, we must verify that the final charge of the bank is nonnegative.
We use the following three discharging rules:
\begin{itemize}
\item[(R1)] Each $\Delta$-vertex and $(\Delta-1)$-vertex $v$ gives a charge of $1/3$ to each adjacent 2-vertex or 3-vertex.
\item[(R2)] Each large vertex $v$ gives a charge of $1/2$ to each incident triangle.
\item[(R3)] Each $\Delta$-vertex gives a charge of $4/3$ to the bank. \\
            Each 2-vertex takes a charge of $4/3$ from the bank.
\end{itemize}
The only rule that effects the bank's charge is (R3).  
Since $G$ has more $\Delta$-vertices than 2-vertices, the bank's final charge is positve.

Now we show that for every element the new charge $\mu^*$ is nonnegative.

\noindent Consider an arbitrary face $f$.  
\begin{itemize}
\item If $d(f) = 3$, then since at least one endpoint of each edge is large, at least two of the vertices incident to $f$ are large.
Thus $\mu^*(f)\geq -1 + 2(1/2) = 0$.  
\item If $d(f)\geq 4$, then $\mu^*(f) = \mu(f)\geq 0$.  
\end{itemize}
Consider an arbitrary vertex $v$.
\begin{itemize}
\item If $d(v) = 2$, then $\mu^*(v) = -2 + 2(1/3) + 4/3 = 0$.  
\item If $d(v) = 3$, then $\mu^*(v) = -1 + 3(1/3) = 0$.  
\item If $d(v) = 4$, then $\mu^*(v) = \mu(v) = 0$.  
\item If $d(v) = 5$, then $v$ is incident to at most 2 triangles, so $\mu^*(v) \geq 1 - 2(1/2) = 0$.  
\item If $6\leq d(v) \leq \Delta-2$, then $v$ is incident to at most $\frac{d(v)}{2}$ triangles.  
Thus $\mu^*(v)\geq {d(v) -4 - \frac{d(v)}{2}(\frac{1}{2}}) = \frac{3d(v)}{4} - 4 > 0$.  
\item If $d(v) = \Delta-1$, then let $t$ be the number of triangles incident to $v$. 
For each triangle incident to $v$, at most one of the vertices of that triangle has degree 3.
Thus, if $v$ is incident to $t$ triangles, then $\mu^*(v)\geq {d(v)-4 -t(\frac{1}{2}) - (d(v) - t)(\frac{1}{3}}) = {\frac{2d(v)}{3} -4 - \frac{t}{6}}$.  
Since $t\leq \frac{d(v)}{2}$, we get $\mu^*(v) \geq \frac{7}{12}d(v) - 4$.  This expression is positive when $d(v)\geq 8$.

\item If $d(v) = \Delta$, then let $t$ be the number of triangles incident to $v$. 
For each triangle incident to $v$, at most one of the vertices of that triangle has degree 3.
Thus, if $v$ is incident to $t$ triangles, then $\mu^*(v)\geq {d(v)-4 -\frac{4}{3} -t(\frac{1}{2}) - (d(v) - t)(\frac{1}{3}}) = {\frac{2d(v)}{3} - \frac{16}{3}- \frac{t}{6}}$.  
Since $t\leq \floor{\frac{d(v)}{2}}$, this expression is nonnegative when $d(v)\geq 9$.
\end{itemize}
\aftermath
\end{proof}

\section{Other Surfaces with Nonnegative Euler Characteristic}
Each of the proofs in Section~2 only used planarity to show that the sum of the initial charges is less than 0.  Thus, each of the proofs also holds for the projective plane (for which $|F(G)|-|E(G)|+|V(G)|=1$).  The proofs also ``almost'' hold for the torus and the Klein bottle (for both surfaces $|F(G)|-|E(G)|+|V(G)|=0$).  Thus, to complete the proofs for the torus and the Klein bottle, it is sufficient to show that after the discharging at least one element has strictly positive charge.  In Theorem~\ref{delta7}, each vertex of degree $\Delta$ has positive charge.  In Theorem~\ref{delta9}, the bank has positive charge.  In Lemma~\ref{delta5-6}, when $\Delta=6$, each vertex of degree $\Delta$ has positive charge.  Since we cannot prove Lemma~\ref{delta5-6} for the 
torus and the Klein bottle when $\Delta=5$, we prove a weaker result:

\begin{lemma}
Let $G$ be a graph with no kites that is embedded in a surface with nonnegative Euler characteristic.  If $\Delta\leq 6$, then $G$ has an edge $uv$ with $d(u)+d(v)\leq 9$.
\end{lemma}
\begin{proof}
If $\Delta\leq 4$, the result holds, since $4+4< 9$.
If $\Delta\in\{5,6\}$, we use the same discharging argument as in Lemma~6.
We have already shown that $\mu^*(x)\geq 0$ for each element $x$; now we
must show there exists some element $x$ such that $\mu^*(x)>0$.

If $\Delta=6$, then each vertex $v$ of degree 6 has $\mu^*(v)>0$.
Suppose $\Delta=5$.  If the lemma fails, then $G$ must be 5-regular.  If $G$ contains no triangles, then each 5-vertex has positive charge; if $G$ contains a triangle $f$, then since $G$ is 5-regular, $\mu^*(f)=-1 + 3(1/2) = 1/2 > 0$.
\end{proof}

We now have the necessary tools to prove our main results. 
In the next section we prove two theorems about edge-choosability.
In the following section we prove two theorems about total-choosability.
For simplicity, we state these theorems only for planar graphs, but after each proof we note which results hold for other surfaces.

\section{Application to Edge-Choosability}

\begin{theorem}
Let $G$ be a planar graph that contains no kites.  
If $\Delta\neq5$, then $\chi'_l(G)\leq\Delta+1$.
If $\Delta=5$, then $\chi'_l(G)\leq\Delta+2$. 
\label{main1}
\end{theorem}
\begin{proof}
Let $G$ be a connected graph.
Harris \cite{harris} and Juvan et al. \cite{juvan1} showed that $G$ is $(\Delta+1)$-edge-choosable when $\Delta=3$ and $\Delta=4$, respectively (even for nonplanar graphs).  
Thus, we only need to prove the theorem when $\Delta\geq 5$.  
We consider separately the three cases $\Delta = 5$, $\Delta = 6$, and $\Delta\geq 7$.  
In each case we proceed by induction on the number of edges.  
The theorem holds trivially if $|E(G)|\leq 7$.
Note that if $d(u)+d(v)\leq k$, then edge $uv$ is adjacent to at most $k-2$ other edges.  We use this fact frequently in the proof.

Suppose $\Delta(G) = 5$.  
Let $H$ be a subgraph of $G$. Since $\Delta(H)\leq 5$, Lemma \ref{delta5-6} implies that $H$ has an edge $uv$ with $d(u) + d(v)\leq 8$.  
By hypothesis, $\chi'_l(H-uv)\leq 7$.  Since edge $uv$ is adjacent to at most six edges in $H$, we can extend the coloring to edge $uv$.

Suppose $\Delta(G) \geq 7$.  
Let $H$ be a subgraph of $G$.  Since $\Delta(H)\leq \Delta(G)$, Theorem~\ref{delta7} and Lemma~\ref{delta5-6} together imply that $H$ has an edge $uv$ with $d(u) + d(v)\leq \Delta(G)+2$.  
By hypothesis, $\chi'_l(H-uv)\leq \Delta(G)+1$.  
Since edge $uv$ is adjacent to at most $\Delta(G)$ edges in $H$, we can extend the coloring to edge $uv$.

Suppose $\Delta(G) = 6$.  
Let $H$ be a subgraph of $G$.  By Lemmas~\ref{delta6} and \ref{delta5-6}, we know that one of the three conditions from Lemma~\ref{delta6} holds for $H$.
We show that in each case 
we can remove some set of edges $\widehat{E}$, inductively color the graph $H-\widehat{E}$, then extend the coloring to $\widehat{E}$.
\smallskip

(i) If $H$ has an edge $uv$ with $d(u)+d(v)\leq 8$, then by hypothesis $\chi'_l(H-uv)\leq 7$.  
Since at most 6 colors are prohibited from use on $uv$, we can extend the coloring to $uv$. 
\smallskip

(ii) If $H$ has a 4-face $uvwx$ with $d(u) = d(w) = 3$, then let $\mathcal{C}=\{uv,vw,wx,xu\}$.  
By hypothesis $\chi'_l(H-\mathcal{C})\leq 7$.  
Since each of the four uncolored edges of $\mathcal{C}$ has at most 5 colors prohibited, there are at least two colors available to use on each edge of $\mathcal{C}$.  
Since $\chi'_l(\mathcal{C}) = 2$, we can extend the coloring to $\mathcal{C}$. 
(It is well-known for every even cycle $\mathcal{C}$ that $\chi'_l(\mathcal{C})=2$; for example, this was shown by Erd\H{o}s, Rubin, and Taylor~\cite{ERT}.)

\smallskip

(iii) If $G$ has a 6-vertex incident to 3 triangles, two of type $(6,6,3)$ and the third of type $(6,6,3)$, $(6,5,4)$, or $(6,6,4)$, then we assume the third triangle is type $(6,6,4)$, since this is the most restrictive case.  
Let $\hat{E}$ be the set of edges of all three triangles, plus one additional edge incident to a vertex of degree 3 in one of the triangles.  
By hypothesis, $\chi'_l(G-\widehat{E})\leq 7$.  We show that we can extend the coloring to $\widehat{E}$.

\begin{figure}
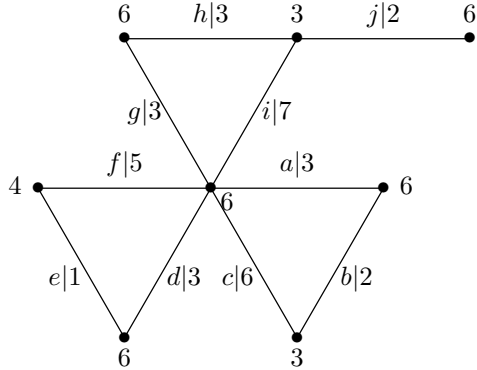

\gpic{
\expandafter\ifx\csname graph\endcsname\relax \csname newbox\endcsname\graph\fi
\expandafter\ifx\csname graphtemp\endcsname\relax \csname newdimen\endcsname\graphtemp\fi
\setbox\graph=\vtop{\vskip 0pt\hbox{%
    \graphtemp=.5ex\advance\graphtemp by 0.120in
    \rlap{\kern 1.476in\lower\graphtemp\hbox to 0pt{\hss $\bullet$\hss}}%
    \graphtemp=.5ex\advance\graphtemp by 0.903in
    \rlap{\kern 1.928in\lower\graphtemp\hbox to 0pt{\hss $\bullet$\hss}}%
    \graphtemp=.5ex\advance\graphtemp by 1.686in
    \rlap{\kern 1.476in\lower\graphtemp\hbox to 0pt{\hss $\bullet$\hss}}%
    \graphtemp=.5ex\advance\graphtemp by 1.686in
    \rlap{\kern 0.572in\lower\graphtemp\hbox to 0pt{\hss $\bullet$\hss}}%
    \graphtemp=.5ex\advance\graphtemp by 0.903in
    \rlap{\kern 0.120in\lower\graphtemp\hbox to 0pt{\hss $\bullet$\hss}}%
    \graphtemp=.5ex\advance\graphtemp by 0.120in
    \rlap{\kern 0.572in\lower\graphtemp\hbox to 0pt{\hss $\bullet$\hss}}%
    \special{pn 8}%
    \special{pa 1476 120}%
    \special{pa 572 1686}%
    \special{pa 120 903}%
    \special{pa 1928 903}%
    \special{pa 1476 1686}%
    \special{pa 572 120}%
    \special{pa 1476 120}%
    \special{fp}%
    \special{pa 1476 120}%
    \special{pa 2380 120}%
    \special{fp}%
    \graphtemp=.5ex\advance\graphtemp by 0.120in
    \rlap{\kern 2.380in\lower\graphtemp\hbox to 0pt{\hss $\bullet$\hss}}%
    \graphtemp=.5ex\advance\graphtemp by 0.903in
    \rlap{\kern 1.024in\lower\graphtemp\hbox to 0pt{\hss $\bullet$\hss}}%
    \graphtemp=.5ex\advance\graphtemp by 0.783in
    \rlap{\kern 1.476in\lower\graphtemp\hbox to 0pt{\hss $a|3$\hss}}%
    \graphtemp=.5ex\advance\graphtemp by 1.379in
    \rlap{\kern 1.787in\lower\graphtemp\hbox to 0pt{\hss $b|2$\hss}}%
    \graphtemp=.5ex\advance\graphtemp by 1.379in
    \rlap{\kern 1.165in\lower\graphtemp\hbox to 0pt{\hss $c|6$\hss}}%
    \graphtemp=.5ex\advance\graphtemp by 1.379in
    \rlap{\kern 0.883in\lower\graphtemp\hbox to 0pt{\hss $d|3$\hss}}%
    \graphtemp=.5ex\advance\graphtemp by 1.379in
    \rlap{\kern 0.261in\lower\graphtemp\hbox to 0pt{\hss $e|1$\hss}}%
    \graphtemp=.5ex\advance\graphtemp by 0.783in
    \rlap{\kern 0.572in\lower\graphtemp\hbox to 0pt{\hss $f|5$\hss}}%
    \graphtemp=.5ex\advance\graphtemp by 0.512in
    \rlap{\kern 0.678in\lower\graphtemp\hbox to 0pt{\hss $g|3$\hss}}%
    \graphtemp=.5ex\advance\graphtemp by 0.000in
    \rlap{\kern 1.024in\lower\graphtemp\hbox to 0pt{\hss $h|3$\hss}}%
    \graphtemp=.5ex\advance\graphtemp by 0.512in
    \rlap{\kern 1.370in\lower\graphtemp\hbox to 0pt{\hss $i|7$\hss}}%
    \graphtemp=.5ex\advance\graphtemp by 0.000in
    \rlap{\kern 1.928in\lower\graphtemp\hbox to 0pt{\hss $j|2$\hss}}%
    \graphtemp=.5ex\advance\graphtemp by 0.903in
    \rlap{\kern 2.048in\lower\graphtemp\hbox to 0pt{\hss 6\hss}}%
    \graphtemp=.5ex\advance\graphtemp by 1.806in
    \rlap{\kern 1.476in\lower\graphtemp\hbox to 0pt{\hss 3\hss}}%
    \graphtemp=.5ex\advance\graphtemp by 1.806in
    \rlap{\kern 0.572in\lower\graphtemp\hbox to 0pt{\hss 6\hss}}%
    \graphtemp=.5ex\advance\graphtemp by 0.903in
    \rlap{\kern 0.000in\lower\graphtemp\hbox to 0pt{\hss 4\hss}}%
    \graphtemp=.5ex\advance\graphtemp by 0.000in
    \rlap{\kern 0.572in\lower\graphtemp\hbox to 0pt{\hss 6\hss}}%
    \graphtemp=.5ex\advance\graphtemp by 0.000in
    \rlap{\kern 1.476in\lower\graphtemp\hbox to 0pt{\hss 3\hss}}%
    \graphtemp=.5ex\advance\graphtemp by 0.000in
    \rlap{\kern 2.380in\lower\graphtemp\hbox to 0pt{\hss 6\hss}}%
    \graphtemp=.5ex\advance\graphtemp by 0.988in
    \rlap{\kern 1.109in\lower\graphtemp\hbox to 0pt{\hss 6\hss}}%
    \hbox{\vrule depth1.806in width0pt height 0pt}%
    \kern 2.500in
  }%
}%
}
\caption{The ten remaining uncolored edges.  The number at each vertex is the degree of that vertex in $G$.  The number on each edge is the number of colors available to use on that edge after we have chosen colors for all edges not pictured.}
\label{fig:ten-edges}
\end{figure}
%
The ten edges of $\widehat{E}$ are shown in Figure~\ref{fig:ten-edges}, along with the number of colors available to use on each edge.  
We use $L(e)$ to denote the list of colors available for use on edge $e$ after we have chosen colors for all the edges not shown in Figure~\ref{fig:ten-edges}.
%
Since $|L(g)| + |L(j)| > |L(h)|$, either there exists some color $\alpha\in L(g) \cap L(j)$ or there exists some color $\alpha\in (L(g)\cup L(j))\setminus L(h)$.  If $\alpha\in L(g)\cap L(j)$, we use color $\alpha$ on edges $g$ and $j$.  Otherwise there exists $\alpha\in (\L(g)\cup L(j))\setminus L(h)$.
In this case, use color $\alpha$ on $g$ or $j$, then use some other available color on whichever of $g$ and $j$ is uncolored.  
In either case, we can now color the rest of the edges in the order: $e,d,a,b,f,c,i,h$.

This completes the proof for the case $\Delta(G) = 6$.
\end{proof}

The results in Theorem~\ref{main1} easily extend to the projective plane for all values of $\Delta$; they also extend to the torus and Klein bottle when $\Delta\geq 7$.

\begin{theorem}
\label{main2}
If $G$ is planar, $G$ contains no kites, and $\Delta(G)\geq 9$, then $\chi'_l(G) = \Delta(G)$.
\end{theorem}
\begin{proof}
Since edges with a common endpoint must receive distinct colors, 
$\chi'_l(G)\geq \Delta(G)$.
So we need to prove that $\chi'_l(G)\leq \Delta(G)$.  
By induction on the number of edges, we prove that if $H$ is a subgraph of $G$, then $\chi'_l(H)\leq\Delta(G)$.  
Our base case is when $\Delta(H)\leq 8$.
The result holds for the base case by Theorem~\ref{main1}.

Assume that $\Delta(H) \geq 9$.
By Theorem~\ref{delta9} at
least one of the following two conditions holds:
\begin{enumerate}
\item[(i)] $H$ has an edge $uv$ with $d(u)+d(v)\leq \Delta(H)+1$.
\item[(ii)] $H$ has a 2-alternating cycle.
\end{enumerate}

Suppose condition (i) holds.  By hypothesis, $\chi'_l(H-uv)\leq \Delta(G)$.  Since $d(u)+d(v)\leq\Delta(H)+1\leq\Delta(G)+1$,
we have at least one color available to extend the coloring to $uv$.

Suppose condition (ii) holds.  Let $\mathcal{C}$ be the even cycle.  By hypothesis, $\chi'_l(H-\mathcal{C})\leq\Delta(G)$.
After coloring $H-\mathcal{C}$, each edge of $\mathcal{C}$ has at least two colors available.  Since even cycles are 
2-choosable, we can extend the coloring to $\mathcal{C}$.
\end{proof}

Theorem~\ref{main2} easily extends to the projective plane, Klein bottle, and torus.

\section{Application to Total-Choosability}
\begin{theorem}
\label{total-color-low}
If $G$ is a planar graph with no kites and $\Delta\geq 7$, then $\chi''_l(G)\leq\Delta+2$.
\end{theorem}
\begin{proof}
Our proof is by induction on the number of edges in $G$.
(Our base case is all planar graphs with at most 4 edges.)
By Theorem~4 and Lemma~8, every subgraph of $G$ contains an edge $uv$ with $d(u)\leq 4$ and $d(u)+d(v)\leq \Delta+2$.
By induction, we have a total coloring for $G-uv$.  Now we uncolor $u$; since edge $uv$ is adjacent to at most $\Delta$ colored edges and one colored vertex, we can color $uv$.  Since $d(u)\leq 4$, $u$ is incident and adjacent to at most $2d(u)\leq 8$ colored elements; thus, we can color $u$.
\end{proof}

Theorem~\ref{total-color-low} extends to the projective plane, Klein bottle, and torus.

\begin{theorem}
If $G$ is a planar graph with no kites and $\Delta\geq 9$, then $\chi''_l(G)=\Delta+1$.
\label{thm15}
\end{theorem}
\begin{proof}
Our proof is by induction on the number of edges in $G$.  Our base case is all subgraphs $H$ of $G$ with $\Delta(H)\in\{7,8\}$.  By Theorem~\ref{total-color-low}, each subgraph $H$ satisfies $\chi''_l(H)\leq\Delta(H)+2\leq 10\leq\Delta(G)+1$.

Our induction step is as follows.
By Lemma~\ref{delta9}, every subgraph $H$ of $G$ with $\Delta(H)\geq 9$ satisfies one of the following two conditions:
(i) $H$ has an edge $uv$ with $d(u)\leq 4$ and $d(u)+d(v)\leq \Delta+1$ or
(ii) $H$ has a 2-alternating cycle. 

Suppose that $H$ satisfies (i).  By the induction hypothesis, we have a total coloring of $H-uv$.  To extend this coloring to $uv$, uncolor vertex $u$; since edge $uv$ is adjacent to at most $\Delta(H)-1$ colored edges and one colored vertex, we can extend the coloring to $uv$.  Finally, since vertex $u$ is adjacent and incident to at most $2d_H(u)\leq 8$ colored elements, we can color $u$.

Suppose instead that $H$ satisfies (ii).  Let $C$ denote the 2-alternating cycle $v_1w_1v_2w_2\ldots v_kw_k$.  By the induction hypothesis, we have a total coloring of $H-E(C)$.  To extend the coloring to $E(C)$, we first uncolor each 2-vertex $w_i$ on $C$.  Now each edge in $E(C)$ has at least two colors available.  Since even cycles are 2-choosable (and equivalently 2-(edge-choosable)), we can color the edges of $C$.  Finally, we can color each 2-vertex $w_i$, since each such vertex is incident and adjacent to only 4 colored elements.
\end{proof}

Theorem~\ref{thm15} easily extends to the projective plane, Klein bottle, and torus.

\section{Acknowledgements} 
Doug West clarified the arguments and improved the exposition.
David Bunde offered many suggestions that improved the exposition.
Oleg Borodin taught a lecture course on discharging at UIUC during the Spring 2005 semester.  Without that course, this paper would not have been possible.
Thanks to an anonymous referee, who provided two very useful references.
Thanks most of all to my Lord and Savior, Jesus Christ.

\end{document}